\documentclass[11pt]{article}

\usepackage[T1]{fontenc}
\usepackage[utf8]{inputenc}
\usepackage[british]{babel}
\usepackage{amsmath,amssymb,amsthm,mathtools}
\usepackage[margin=25mm,a4paper]{geometry}
\usepackage{microtype}
\usepackage{enumitem}
\usepackage{xcolor}
\usepackage[colorlinks=true,linkcolor=blue!55!black,citecolor=blue!55!black,
  urlcolor=blue!55!black]{hyperref}

\newcommand{\ceil}[1]{\left\lceil #1 \right\rceil}
\newcommand{\floor}[1]{\left\lfloor #1 \right\rfloor}
\newcommand{\abs}[1]{\lvert #1\rvert}
\newcommand{\set}[1]{\{#1\}}
\newcommand{\nonN}{\overline N}

\theoremstyle{plain}
\newtheorem{maintheorem}{Theorem}
\newtheorem{theorem}{Theorem}[section]

\newtheorem{lemma}[theorem]{Lemma}

\theoremstyle{definition}

\newtheorem{conjecture}{Conjecture}
\newtheorem{problem}{Problem}

\theoremstyle{remark}

\numberwithin{equation}{section}
\allowdisplaybreaks

\title{Totally odd clique immersions in graphs with bounded independence number}
\author{Sebasti\'an Bustamante\thanks{Independent researcher. Email: \href{mailto:dev@sebustam.com}{dev@sebustam.com}.}}
\date{}

\begin{document}
\maketitle

\begin{abstract}
We prove that every $n$-vertex graph with at least one edge and independence
number at most an integer $a\ge2$, where $n\ge a$, contains a totally odd
strong immersion of
\[
 K_{\left\lceil (n+a-1)/(2a-3/2) \right\rceil},
\]
with every branch path of length at most $2a-1$. We also prove that every
graph $G$ with independence number at most $a$ and no totally odd strong
$K_t$-immersion with branch paths of length at most $2a-1$ satisfies
\[
 \chi(G)\le\left\lfloor \frac{(2a-3/2)(t-1)-a+2}{2} \right\rfloor,
\]
where $t\ge2$.
Under the same hypotheses, we also obtain
$\chi(G)\le3\log a\,(t-1)$, where $\log$ is the natural logarithm.
The order bound narrows the gap to the conjectured clique order
$\left\lceil n/\alpha(G) \right\rceil$, reaching it up to a factor of 2.
For the conjectured clique order $\chi(G)$, the first colouring bound
gives the sharper guarantee for small $a$, while the second colouring
bound reaches $\chi(G)$ up to a logarithmic factor in $a$.
\end{abstract}

\section{Introduction}
\label{sec:introduction}

All graphs are finite and simple. For a graph $G$, write $n=\abs{V(G)}$
and let $\alpha(G)$ denote its independence number, with $\alpha(G)=0$
for the empty graph. We write
$[a]=\set{1,\ldots,a}$ and
\[
 \nonN_G(v)=V(G)\setminus\bigl(N_G(v)\cup\set v\bigr).
\]
Subscripts are omitted when the host graph is unambiguous.

A \emph{weak immersion} of $H$ in $G$ assigns distinct \emph{branch vertices}
of $G$ to the vertices of $H$ and pairwise edge-disjoint \emph{branch paths}
to its edges, joining the corresponding branch vertices.  It is
\emph{strong} if no branch vertex is internal to a branch path, and it is
\emph{totally odd} if every branch path has odd length.  Unless stated
otherwise, immersions are strong.

For nonempty $G$, let $\iota_o(G)$ be the largest order of a totally
odd strong clique immersion in $G$, and let $\iota_o^{\le L}(G)$ be the
largest such order in which every branch path has length at most $L$.
These parameters are monotone under taking subgraphs, and
$\iota_o^{\le L}(G)\le\iota_o(G)$.

Clique-minor bounds have been studied in terms of order and independence
number~\cite{DM1982,BK2011} and in terms of chromatic
number~\cite{OpenAI2026Hadwiger}.
For immersions, Vergara~\cite{Vergara2017} conjectured a weak
$K_{\ceil{n/2}}$-immersion when $\alpha(G)=2$.
This follows from Deng and Wu~\cite{DW2026}, who proved that every graph
$G$ with $\alpha(G)\le2$ has a strong $K_{\chi(G)}$-immersion.
More generally,
Bustamante, Quiroz, Stein, and Zamora~\cite{BQSZ2022} formulated the
following conjecture.

\begin{conjecture}[BQSZ~\cite{BQSZ2022}]
\label{conj:bqsz}
Every nonempty $n$-vertex graph $G$ contains a weak immersion of
$K_{\ceil{n/\alpha(G)}}$.
\end{conjecture}

They proved in~\cite{BQSZ2022} that every $n$-vertex graph $G$ with
independence number $a\ge3$ contains a totally odd strong clique
immersion of order
\[
 \frac{n}{(9/4)a-O(\log a)}-O(1).
\]
The following theorem reduces the leading denominator coefficient from
$9/4$ to $2$.

\begin{maintheorem}[Order bound]
\label{thm:immersion-order}
Let $a\ge2$ be an integer. Every $n$-vertex graph $G$ with at least
one edge and $\alpha(G)\le a$ contains a totally odd strong immersion of
\[
 K_{\ceil{(n+a-1)/(2a-3/2)}},
\]
when $n\ge a$; when $n<a$, the edge itself gives a $K_2$-immersion.
Every branch path can be chosen to have length at most $2a-1$.
\end{maintheorem}

Theorem~\ref{thm:immersion-order} narrows the gap to the clique order
$\ceil{n/\alpha(G)}$ sought in Conjecture~\ref{conj:bqsz}, while
guaranteeing a totally odd strong immersion with branch paths of length
at most $2a-1$.
Without the total-oddness requirement, using the bound of Deng and
Wu~\cite{DW2026} in the routing argument reduces the denominator
from $2a-3/2$ to $2a-2$.

For independence number two, Vergara~\cite{Vergara2017} proved a
strong clique immersion of order $\ceil{n/3}$. Gauthier, Le, and
Wollan~\cite{GLW2019} obtained a totally odd strong clique immersion of
order $2\floor{n/5}$.

The chromatic question asks whether every graph $G$ contains an
immersion of $K_{\chi(G)}$. Lescure and Meyniel~\cite{LM89} conjectured
the strong version, and Abu-Khzam and Langston~\cite{AKL03} proposed
the weak version.
Churchley~\cite{Churchley2017} proposed a version using pairwise
edge-disjoint odd trails, and Jim\'enez, Quiroz, and Thraves
Caro~\cite{JQTC2024} asked for the totally odd strong version using paths:

\begin{conjecture}[JQTC~\cite{JQTC2024}]
\label{conj:jqtc}
Every graph $G$ contains a totally odd strong immersion of $K_{\chi(G)}$.
\end{conjecture}

Echeverr\'ia and McDonald~\cite{EM2026} proved that if $\alpha(G)\le2$
and $G$ has no totally odd strong $K_t$-immersion, then
$\chi(G)\le\ceil{3(t-1)/2}$.
The following theorem gives a colouring bound for arbitrary $a\ge2$
while controlling the lengths of the branch paths.

\begin{maintheorem}[Linear colouring bound]
\label{thm:colouring-introduction}
Let $a,t\ge2$ be integers. If $\alpha(G)\le a$ and $G$ has no totally odd
strong $K_t$-immersion whose branch paths have length at most $2a-1$,
then
\[
 \chi(G)\le\floor{\frac{(2a-3/2)(t-1)-a+2}{2}}.
\]
\end{maintheorem}

For $\alpha(G)=2$, Theorem~\ref{thm:colouring-introduction} reduces the
coefficient of $t-1$ in the colouring bound in~\cite{EM2026} from $3/2$
to $5/4$.

The next theorem improves the asymptotic dependence on $a$ from linear
in Theorem~\ref{thm:colouring-introduction} to logarithmic.

\begin{maintheorem}[Logarithmic colouring bound]
\label{thm:log-colouring}
Let $a,t\ge2$ be integers. If $\alpha(G)\le a$ and $G$ has no totally odd
strong $K_t$-immersion whose branch paths have length at most $2a-1$,
then
\[
 \chi(G)\le3\log a\,(t-1).
\]
\end{maintheorem}

Theorem~\ref{thm:log-colouring} guarantees the clique order $\chi(G)$
sought in Conjecture~\ref{conj:jqtc} up to a logarithmic factor in $a$:
every nonempty graph $G$ with $\alpha(G)\le a$ contains a totally odd
strong immersion of $K_{\ceil{\chi(G)/(3\log a)}}$ with branch paths
of length at most $2a-1$.

The routing lemmas in Section~\ref{sec:routing} let us extend immersion
guarantees from smaller independence numbers, yielding
Theorem~\ref{thm:immersion-order} in Section~\ref{sec:large-immersions}.
The resulting upper bound on the order of graphs excluding the prescribed
immersion gives the two colouring bounds: Section~\ref{sec:colouring}
uses Gallai's theorem, while Section~\ref{sec:log-colouring} successively
removes maximum independent sets. Appendix~\ref{sec:two-fifths} proves
the order bound for $\alpha(G)=2$.

\section{Routing lemmas}
\label{sec:routing}

We use the following vertex-set form of Menger's theorem~\cite{Menger1927}:
for disjoint vertex sets $S,T$ in a finite graph, the maximum number of
pairwise vertex-disjoint $S$--$T$ paths equals the minimum size of a set
$C$ such that deleting $C$ leaves no path from $S\setminus C$ to
$T\setminus C$. In particular, the separator may contain vertices of
$S\cup T$.

Let $I=\set{a_1,\ldots,a_a}$ be independent, and let $Q,M$ be disjoint
sets outside $I$.  For $J\subseteq[a]$, write
\[
 N(J)=\bigcup_{i\in J}N(a_i),
 \qquad h(J)=\abs{M\setminus N(J)}.
\]
For a partition $[a]=J\mathbin{\dot\cup}J'$, put
\[
 Q(J,J')=Q\cap N(J)\cap N(J').
\]
The paths in the next lemma use only edges with one endpoint in $I$
and the other in $Q\cup M$.  Their disjointness condition allows them
to share vertices of $I$.

\begin{lemma}
\label{lem:terminal-cut}
Suppose every vertex of $M$ has a neighbour in $I$, and fix
$a_\ell\in I$.  There are paths from $a_\ell$ to the distinct vertices
of $M$, alternating between $I$ and $Q\cup M$, pairwise vertex-disjoint
outside $I$, and with all internal vertices in
$Q\cup(I\setminus\set{a_\ell})$, if and only if
\begin{equation}
\label{eq:terminal-cut}
  \abs{Q(J,J')}\ge h(J)
\end{equation}
for every partition $[a]=J\mathbin{\dot\cup}J'$ with
$\ell\in J$ and $J'\ne\varnothing$.  These paths are edge-disjoint
and have odd length at most $2a-1$.
\end{lemma}

\begin{proof}
Suppose first that the required paths exist. We show that
$\abs{Q(J,J')}\ge h(J)$ for every partition
$[a]=J\mathbin{\dot\cup}J'$ with $\ell\in J$ and $J'\ne\varnothing$.
Each path to a terminal in $M\setminus N(J)$ must cross from an
$I$-vertex indexed by $J$ to one indexed by $J'$, through a vertex
of $Q(J,J')$. Distinct paths use distinct such vertices, proving
the inequality.

Conversely, suppose that~\eqref{eq:terminal-cut} holds for every such partition.
We construct the required paths.
Put $r=\abs M$; the case $r=0$ is immediate.
Replace each $a_i$ by an independent bag $A_i$ of $r$ twins in a
bipartite graph with other side $Q\cup M$, retaining precisely the
adjacencies prescribed by $G$. We use the vertex-set form of
Menger's theorem~\cite{Menger1927}, in which a separator may contain
vertices of either endpoint set.
Suppose an $A_\ell$--$M$ separator $C$ has size less than $r$.
Each bag has a surviving vertex. Let $J$ index the bags reachable
from $A_\ell\setminus C$ after deleting $C$.
Here $\ell\in J$, and all surviving vertices of each reachable bag
are reachable, since they are twins.
If $J=[a]$, every surviving terminal is reachable, forcing
$M\subseteq C$, a contradiction.
Otherwise, writing $J'=[a]\setminus J$, the separator contains both
$Q(J,J')$ and $M\cap N(J)$. Thus
\[
 \abs C\ge\abs{Q(J,J')}+\abs{M\cap N(J)}\ge r,
\]
again a contradiction.

Menger's theorem gives $r$ vertex-disjoint paths. Trim them at their
first vertices of $M$, and identify the twins in each bag with $a_i$.
Any repeated vertex of a resulting walk lies in $I$; deleting the
segment between two occurrences preserves alternation and oddness.
The resulting simple paths remain disjoint outside $I$, so they are
edge-disjoint. Their internal vertices avoid $M$ and the root.
Each uses at most $a$ vertices of $I$ and has length at most $2a-1$.
\end{proof}

To apply Lemma~\ref{lem:terminal-cut}, we need to show that
$\abs{Q(J,J')}-h(J)$ is nonnegative for every relevant partition.
The next lemma expresses this difference in terms of the sizes
of common non-neighbourhoods.  When $I$ is a maximum independent
set, these common non-neighbourhoods have smaller independence
number.  In the applications below, this allows us to bound their
sizes by induction under the immersion-exclusion hypothesis, and
hence verify the routing condition.

\begin{lemma}
\label{lem:cut-identity}
Suppose $I=\set{a_1,\ldots,a_a}$ is a maximum independent set and
$V(G)=I\mathbin{\dot\cup}M\mathbin{\dot\cup}Q$, with $\abs M=m$.
For nonempty $S\subsetneq[a]$, put
\[
 X_S=\bigcap_{i\in S}\nonN(a_i).
\]
Then $\alpha(G[X_S])\le a-\abs S$.  For every partition
$[a]=J\mathbin{\dot\cup}J'$ with both sides nonempty,
\begin{equation}
\label{eq:cut-identity}
 \abs{Q(J,J')}-h(J)
 =n-\abs{X_J}-\abs{X_{J'}}-m+h(J').
\end{equation}
\end{lemma}

\begin{proof}
The set $X_S$ avoids $\set{a_i:i\in S}$, and every one of its vertices
is nonadjacent to that set.  An independent set in $G[X_S]$ can therefore
be extended by all $\abs S$ of those vertices, proving the independence
bound.

Every vertex outside $I$ has a neighbour in $I$.  Thus $X_J$ and
$X_{J'}$ are disjoint.  They contain all of $I$ between them: a vertex
$a_i$ indexed by $J$ belongs to $X_{J'}$, and conversely.  Outside
$I$, a vertex not in their union has neighbours on both sides.  Put
$M(J,J')=M\cap N(J)\cap N(J')$.  We obtain
\begin{align*}
 n&=\abs{X_J}+\abs{X_{J'}}+\abs{Q(J,J')}+\abs{M(J,J')},\\
 m&=h(J)+h(J')+\abs{M(J,J')}.
\end{align*}
Subtracting the second identity from the first gives
\eqref{eq:cut-identity}.
\end{proof}

\section{Proof of Theorem~\ref{thm:immersion-order}}
\label{sec:large-immersions}

Our goal is to bound the order of a graph $G$ with $\alpha(G)\le a$ and
no totally odd strong $K_{m+1}$-immersion with branch paths of length at
most $2a-1$. For $a=1$, the bound is $m$, since $G$ is complete. For
$m=1$, the bound is $a$, since $G$ is edgeless. The following recurrence
extends a bound for $a=2$ to all $a\ge3$.

\begin{theorem}
\label{thm:routing-recurrence}
Let $W_a(m)$ be integers indexed by $a,m\ge1$, with $W_1(m)=m$ and
$W_a(1)=a$. Suppose that for every $m\ge2$, every graph with independence
number at most two and no totally odd strong $K_{m+1}$-immersion with
branch paths of length at most three has at most $W_2(m)$ vertices.
If, for $a\ge3$ and $m\ge2$,
\[
 W_a(m)\ge\max\left\{
 W_{a-1}(m),\ W_a(m-1)+a,\
 \max_{1\le s<a}\bigl(W_s(m)+W_{a-s}(m)+m-1\bigr)
 \right\},
\]
then for all $a,m\ge1$, every graph $G$ with $\alpha(G)\le a$ and no
totally odd strong $K_{m+1}$-immersion with branch paths of length at
most $2a-1$ has at most $W_a(m)$ vertices.
\end{theorem}

\begin{proof}
Induct first on $m$ and then on $a$. The cases $m=1$ and $a=1$
were established above, and $a=2$ with $m\ge2$ is the assumed input.
Suppose $a\ge3$, $m\ge2$, and $G$ excludes a totally odd
strong $K_{m+1}$-immersion whose branch paths have length at most $2a-1$.
Assume for a contradiction that $n>W_a(m)$.
The first inequality lets us assume $\alpha(G)=a$. Delete a maximum
independent set $I$ of size $a$. Since $n-a>W_a(m-1)$, induction
gives a totally odd strong $K_m$-immersion in $G-I$ whose branch paths
have length at most $2a-1$. Let $M$ be its branch set and put
$Q=V(G)\setminus(I\cup M)$.
For a partition $[a]=J\mathbin{\dot\cup}J'$ into nonempty parts,
\eqref{eq:cut-identity} and induction give
\[
 |Q(J,J')|-h(J)
 \ge n-W_{|J'|}(m)-W_{|J|}(m)-m+h(J')\ge0.
\]
Here the exclusion hypothesis also holds in each induced subgraph, with
the corresponding shorter path-length restriction. Lemma~\ref{lem:terminal-cut}
joins any root in $I$ to
$M$ by edge-disjoint odd paths of length at most $2a-1$, internally
avoiding the branch vertices. Every new edge meets $I$, so none belongs
to the old immersion. This gives the forbidden $K_{m+1}$-immersion.
\end{proof}

We apply the parameterized recurrence of Theorem~\ref{thm:routing-recurrence}
to the following bound for graphs with independence number at most two.

\begin{theorem}
\label{thm:two-fifths}
Every $n$-vertex graph with at least one edge and independence number at most
two contains a totally odd strong immersion of
\[
 K_{\ceil{(2n+2)/5}}
\]
in which every branch path has length at most three.
\end{theorem}

The proof of Theorem~\ref{thm:two-fifths} is deferred to
Appendix~\ref{sec:two-fifths}.

In exclusion form, Theorem~\ref{thm:two-fifths} gives the following bound:
for a graph with independence number at most two and an
integer $m\ge1$, exclusion of a totally odd strong $K_{m+1}$-immersion
whose branch paths have length at most three gives the integer estimate
\begin{equation}
\label{eq:rounded-exclusion}
 n\le
 \begin{cases}
 2,&m=1,\\
 \floor{5m/2}-1,&m\ge2.
 \end{cases}
\end{equation}
If $m=1$, exclusion of $K_2$ makes the graph edgeless, so $n\le2$. For
$m\ge2$, a graph with an edge and no such $K_{m+1}$-immersion satisfies
$\ceil{(2n+2)/5}\le m$, hence $n\le\floor{5m/2}-1$. An edgeless graph
has at most two vertices and also satisfies this bound. This is the
independence-two input used in the parameterized recurrence.

\begin{lemma}
\label{lem:rounded-exclusion}
If $a,m\ge1$ are integers, $\alpha(G)\le a$, and $G$ excludes a totally odd
strong $K_{m+1}$-immersion whose branch paths have length at most $2a-1$,
then
\[
 n\le V_a(m):=
 \begin{cases}
 a,&m=1,\\
 m,&a=1,\\
 (2a-1)m-(a-1)-\ceil{m/2},&a,m\ge2.
 \end{cases}
\]
\end{lemma}

\begin{proof}
Theorem~\ref{thm:two-fifths}, in the exclusion form~\eqref{eq:rounded-exclusion},
supplies the independence-two input.
We check the inequalities in Theorem~\ref{thm:routing-recurrence}.
For $a\ge3,m\ge2$,
$V_a(m)-V_a(m-1)\ge2a-2\ge a$, including $m=2$ since
$V_a(2)=3a-2$. The bounds increase with $a$. For $s+t=a$,
$s,t\ge1$, and $m\ge2$, direct calculation gives
\[
 V_s(m)+V_t(m)+m-V_a(m)=
 \begin{cases}
 1,&\min\{s,t\}=1,\\
 1-\ceil{m/2},&s,t\ge2.
 \end{cases}
\]
Each value is at most one, so all the recurrence inequalities hold.
Theorem~\ref{thm:routing-recurrence} gives the assertion.
\end{proof}

For $a,m\ge2$, the bound in Lemma~\ref{lem:rounded-exclusion}
can be written as
\begin{equation}
\label{eq:integer-rounding}
 V_a(m)=\floor{\left(2a-\frac32\right)m}-a+1.
\end{equation}
Indeed, $(2a-1)m$ is an integer, so subtracting $\ceil{m/2}$
is the same as subtracting $m/2$ and then taking the floor.
Thus, for $k\ge3$, the order threshold supplied by the lemma for
forcing $K_k$ is
\[
 n\ge\floor{\left(2a-\frac32\right)(k-1)}-a+2,
 \quad\text{equivalently}\quad
 n>\left(2a-\frac32\right)(k-1)-a+1.
\]
Inverting this threshold gives Theorem~\ref{thm:immersion-order}.

\begin{proof}[Proof of Theorem~\ref{thm:immersion-order}]
If $n<a$, an edge supplies the required $K_2$-immersion.
Suppose now that $n\ge a$, and put $C=2a-3/2$ and
$k=\ceil{(n+a-1)/C}$. An edge supplies the conclusion for $k\le2$.
If $k\ge3$, exclusion of $K_k$ with the required path-length restriction
would give, by Lemma~\ref{lem:rounded-exclusion}
and \eqref{eq:integer-rounding},
\[
 n\le\floor{C(k-1)}-(a-1)\le C(k-1)-(a-1)<n,
\]
a contradiction.
\end{proof}

\section{Proof of Theorem~\ref{thm:colouring-introduction}}
\label{sec:colouring}

To obtain a colouring bound linear in $a$, we combine the order bound
from Section~\ref{sec:large-immersions} with Gallai's theorem~\cite{Gallai1963}:
every $k$-vertex-critical graph on at most $2k-2$ vertices is a join of
two nonempty graphs. The join $G_1\vee G_2$ adds all edges between its
two vertex sets, and its chromatic number is additive.
The following lemma abstracts the Gallai-criticality and join strategy
used by Echeverr\'ia and McDonald~\cite{EM2026}, converting an
immersion-exclusion order bound into a colouring bound.

\begin{lemma}
\label{thm:gallai-join}
Fix an integer $a\ge2$ and real numbers $c,d$ with
$d\ge0$ and $c-d\ge1$. Suppose that, for every integer $m\ge2$,
every graph $G$ with $\alpha(G)\le a$ excluding a totally odd
strong $K_{m+1}$-immersion whose branch paths have length at most
$2a-1$ satisfies
\[
 |V(G)|\le2(cm-d)-1.
\]
Then every nonempty graph $G$ with $\alpha(G)\le a$ satisfies
\[
 \chi(G)\le\floor{c\iota_o^{\le2a-1}(G)-d}.
\]
\end{lemma}

\begin{proof}
Take a smallest counterexample, and put
$k=\chi(G)$ and $m=\iota_o^{\le2a-1}(G)$.
If $m=1$, then $G$ is edgeless and
$k=1\le\floor{c-d}$, a contradiction. Thus $m\ge2$.
By the definition of $m$, the exclusion hypothesis applies to $G$.
Minimality makes $G$ vertex-critical: otherwise, deleting a vertex
without decreasing the chromatic number gives a smaller counterexample,
since the immersion parameter cannot increase and $c>0$.
Since the claimed inequality
fails, $k>cm-d$, and hence
\[
 |V(G)|\le2(cm-d)-1<2k-1.
\]
Gallai's theorem therefore makes $G$ a join $G_1\vee G_2$ of two
nonempty graphs. Put $m_i=\iota_o^{\le2a-1}(G_i)$ for $i=1,2$.
Both graphs have independence number at most $a$, so minimality gives
$\chi(G_i)\le cm_i-d$. Combining their immersions with direct edges
between the branch sets gives $m\ge m_1+m_2$. Consequently,
\[
 k=\chi(G_1)+\chi(G_2)
 \le c(m_1+m_2)-2d
 \le cm-2d
 \le cm-d,
\]
a contradiction.
\end{proof}

We now apply Lemma~\ref{thm:gallai-join} to obtain the corresponding
colouring bound.

\begin{theorem}
\label{thm:colouring-bound}
Let $a\ge2$ be an integer. Every nonempty graph $G$ with $\alpha(G)\le a$
satisfies
\[
 \chi(G)\le
 \floor{\frac{(2a-3/2)\iota_o^{\le2a-1}(G)-a+2}{2}}.
\]
\end{theorem}

\begin{proof}
Apply Lemma~\ref{thm:gallai-join} with
\[
 c=a-\frac34,\qquad d=\frac{a-2}{2}.
\]
Here $d\ge0$ and $c-d=a/2+1/4\ge1$.
For every integer $m\ge2$,
Lemma~\ref{lem:rounded-exclusion} supplies the required
exclusion bound:
\[
 |V(G)|\le\left(2a-\frac32\right)m-a+1
          =2(cm-d)-1.
\]
Lemma~\ref{thm:gallai-join} gives the asserted colouring bound.
\end{proof}

\begin{proof}[Proof of Theorem~\ref{thm:colouring-introduction}]
Under the exclusion hypothesis of Theorem~\ref{thm:colouring-introduction},
$\iota_o^{\le2a-1}(G)\le t-1$, so Theorem~\ref{thm:colouring-bound}
gives the exclusion bound for nonempty graphs.
The empty graph satisfies the exclusion conclusion trivially.
\end{proof}

\section{Proof of Theorem~\ref{thm:log-colouring}}
\label{sec:log-colouring}

The bound in Theorem~\ref{thm:colouring-introduction} grows linearly with
$a$ for fixed $t$. To obtain a logarithmic bound in $a$, we successively
remove maximum independent sets and apply Lemma~\ref{lem:rounded-exclusion}
to bound the order of the remaining graph as its independence number
decreases. This gives the following bound in terms of harmonic numbers.

\begin{theorem}
\label{thm:harmonic-colouring}
Let $a,t\ge2$ be integers. If $\alpha(G)\le a$ and $G$ has no totally odd
strong $K_t$-immersion whose branch paths have length at most $2a-1$,
then for $t\ge3$,
\[
 \chi(G)\le
 \floor{\frac{B_1+B_2}{2}+(2t-3)(H_a-H_2)},
\]
where $B_1=t-1$, $B_2=\floor{5(t-1)/2}-1$, and
$H_a=\sum_{j=1}^a 1/j$. If $t=2$, then $\chi(G)\le1$.
\end{theorem}

\begin{proof}
If $t=2$, then $G$ is edgeless, so $\chi(G)\le1$.
Suppose $t\ge3$. For $1\le j\le a$, the order bound from
Lemma~\ref{lem:rounded-exclusion} is
\[
 B_j=V_j(t-1)=
 \begin{cases}
  t-1, & j=1,\\
  \floor{(2j-\tfrac32)(t-1)}-j+1, & 2\le j\le a.
 \end{cases}
\]
Every induced subgraph of $G$ with independence number at most $j$ has
at most $B_j$ vertices by Lemma~\ref{lem:rounded-exclusion}.
Indeed, any totally odd strong $K_t$-immersion with branch paths of length
at most $2j-1$ in the subgraph would also satisfy the length bound $2a-1$
in $G$.

For $j=a,a-1,\ldots,2$, remove independent sets of size $j$ until none
remain, assigning a distinct colour to each set. Let $n_j$ be the order
of the graph at the start of stage $j$, and let $n_1$ be the order of the
final graph, which is a clique or is empty. At the start of stage $j$ the
independence number is at most $j$, so $n_j\le B_j$ for all $j$.
Giving each remaining vertex a separate colour, we use
\[
 n_1+\sum_{j=2}^a\frac{n_j-n_{j-1}}{j}
 =\frac{n_a}{a}+\sum_{j=1}^{a-1}\frac{n_j}{j(j+1)}
\]
colours. Since $n_j\le B_j$ and the coefficients on the right are positive,
the same summation identity gives
\[
 \begin{aligned}
 \chi(G)&\le\frac{B_a}{a}+\sum_{j=1}^{a-1}\frac{B_j}{j(j+1)}\\
 &=\frac{B_1+B_2}{2}+\sum_{j=3}^a\frac{B_j-B_{j-1}}{j}\\
 &=\frac{B_1+B_2}{2}+(2t-3)(H_a-H_2),
 \end{aligned}
\]
where $B_j-B_{j-1}=2t-3$ for $j\ge3$ by Lemma~\ref{lem:rounded-exclusion}.
Taking the floor proves the assertion.
\end{proof}

\begin{proof}[Proof of Theorem~\ref{thm:log-colouring}]
If $t=2$, Theorem~\ref{thm:harmonic-colouring} gives
$\chi(G)\le1\le3\log a$. For $t\ge3$, we have
$(B_1+B_2)/2\le\tfrac74(t-1)$ and
\[
 H_a-H_2\le\int_2^a\frac{dx}{x}=\log(a/2).
\]
By convexity of $1/x$, $\log 2=\int_1^2 dx/x\ge2/3$.
Theorem~\ref{thm:harmonic-colouring} therefore gives
\[
 \chi(G)\le\bigl(\tfrac74+2\log(a/2)\bigr)(t-1)
 \le 3\log a\,(t-1),
\]
as $\tfrac74\le3\log 2$ and $\log(a/2)\ge0$.
\end{proof}

\section{Discussion}

The bounds proved here admit some marginal refinements. Pursuing them
would complicate the proofs without shedding further light on the
methods. The present arguments appear to offer no further improvement
in the asymptotics highlighted in the abstract: the factor of $2$ in
the order bound and the logarithmic dependence on $a$ in the colouring
bound.

For graphs with independence number two, Gauthier, Le, and
Wollan~\cite{GLW2019} guarantee a totally odd strong clique immersion
of order $2\floor{n/5}$. Conjecture~\ref{conj:bqsz} predicts a
coefficient of $1/2$ for weak immersions, and
Conjecture~\ref{conj:jqtc} implies the same coefficient for totally odd
strong immersions, since $\chi(G)\ge\ceil{n/2}$. Two disjoint cliques of
equal order show that $1/2$ is best possible.

\begin{problem}
Find a constant $c\in(2/5,1/2]$ such that, for all sufficiently large
$n$, every $n$-vertex graph $G$ with $\alpha(G)=2$ contains a totally
odd strong clique immersion of order at least $cn$.
\end{problem}

A stronger independence-two bound would improve the constant term in
the denominator $2a-3/2$, reaching $-2$ if $c=1/2$. This would
strengthen the order guarantee as $n$ grows with $a$ fixed, while
leaving the coefficient of $a$ equal to $2$.

For general independence number, reducing the leading coefficient of
$a$ in the denominator $2a-3/2$ would narrow the asymptotic gap to the
conjectured order $\ceil{n/\alpha(G)}$. This appears to require ideas
beyond the present routing argument.

\begin{problem}
Find a constant $c\in[1,2)$ such that every $n$-vertex graph $G$ with
sufficiently large independence number contains a totally odd strong
clique immersion of order at least $n/(c\alpha(G))$.
\end{problem}

A smaller coefficient of $a$ in the denominator of the immersion
order bound would sharpen both colouring bounds, while their dependence
on $a$ would remain linear and logarithmic, respectively.
Conjecture~\ref{conj:jqtc} predicts $\chi(G)\le t-1$, independently of $a$.

\begin{problem}
Find a function $f(a)=o(\log a)$ such that, for all integers
$a,t\ge2$, every graph $G$ with $\alpha(G)\le a$ and no totally odd
strong $K_t$-immersion satisfies $\chi(G)\le f(a)(t-1)$.
\end{problem}

A bound of the form $\chi(G)\le C(t-1)$ for an absolute constant
$C>1$ would remove the dependence on $a$, while falling short of the
conjectured bound $\chi(G)\le t-1$.

\subsection*{Declaration of AI use}
AI tools were used to assist with proof exploration, refinement
and exposition. The author takes full responsibility for the
mathematical content and the final manuscript.

\bibliographystyle{plain}
\bibliography{refs}

\clearpage
\appendix
\section*{Appendix}
\section{Proof of Theorem~\ref{thm:two-fifths}}
\label{sec:two-fifths}

We adapt the five-cycle deletion argument of
Gauthier--Le--Wollan~\cite{GLW2019} to obtain the rounded target
$\ceil{(2n+2)/5}$, while keeping every branch path odd and of
length at most three. This rounding improvement supplies the
additive saving in the immersion-exclusion bound needed to meet
the order hypothesis of Gallai's theorem in
Section~\ref{sec:colouring}.

Deleting an induced five-cycle reduces the target order by two.
We extend an immersion in the remaining graph by two branch
vertices. The improved rounding requires additional base cases:
Subsection~\ref{sec:appendix-small-orders} establishes the
immersions at orders five, seven, nine, and twelve.
Subsection~\ref{sec:appendix-induction} proves the five-cycle
extension lemma and completes the induction.

Throughout this appendix, all graphs have independence number at most
two. Thus the nonneighbours of a vertex form a clique, and every
vertex outside a nonadjacent pair is adjacent to at least one member
of the pair.

\subsection{Small order immersions}
\label{sec:appendix-small-orders}

\begin{lemma}
\label{lem:odd-pair}
Let $x,y$ be nonadjacent, let $C=N(x)\cap N(y)$, and suppose that
$G-\{x,y\}$ has a totally odd strong $K_m$-immersion with branch
set $M$ and branch paths of length at most three.
It extends to a totally odd strong $K_{m+1}$-immersion with the
same path-length bound if
\[
  2\abs C-\abs{C\cap M}\ge m-1.
\]
In particular, it suffices that $\abs C\ge m-1$, or that
$\abs C\ge m-2$ and $C\setminus M\ne\varnothing$.
\end{lemma}

\begin{proof}
Put $r=\abs{M\cap C}$, and let $p,q$ count the branch vertices
adjacent only to $x$, only to $y$, respectively. Then $p+q+r=m$.
If both $p,q$ exceed $\abs C-r$, then
$m\ge2\abs C-r+2$, a contradiction. By symmetry, suppose
$q\le\abs{C\setminus M}$. Add $x$ as a branch vertex.
Use direct edges to its old branch neighbours, and for each other
branch vertex $u$ use $x,z_u,y,u$, with distinct
$z_u\in C\setminus M$. All new paths have length one or three,
internally avoid the new branch set, and are edge-disjoint.
Every new edge meets $\{x,y\}$, so none belongs to the old immersion.
\end{proof}

We use the following observation about a triangle-free
nonbipartite graph $F$. If $Q$ is a shortest odd cycle, then $Q$ is
induced, every vertex outside $Q$ has at most two neighbours on $Q$,
and hence
\begin{equation}\label{eq:cycle-degree}
  \sum_{v\in V(Q)}d_F(v)\le2\abs{V(F)}.
\end{equation}
Indeed, three or more neighbours on $Q$ have cyclic gaps of length
at least two. An odd gap has length at most $\abs{Q}-4$ and gives
a shorter odd cycle through the outside vertex.

\begin{lemma}
\label{lem:odd-small}
Every $n$-vertex graph with independence number at most two has a
totally odd strong $K_s$-immersion with branch paths of length
at most three for
\[
  (n,s)\in\{(5,3),(7,4),(9,5),(12,6)\}.
\]
\end{lemma}

\begin{proof}
We prove the cases in the displayed order; each conclusion also
holds at larger orders by taking an induced subgraph.

\emph{Five vertices.}
The graph is not bipartite. A shortest odd cycle has length three
or five; three consecutive cycle vertices give the required
immersion.

For each subsequent case, suppose there is a counterexample $G$.
Take an edge-minimal spanning subgraph with $\alpha(G)\le2$.
Then $F=\overline G$ is maximal triangle-free. A clique of order
$s$ in $G$ would suffice, so
\[
  \alpha(F)\le s-1,\qquad \Delta(F)\le s-1.
\]
Also, $F$ has no isolated vertex: the corresponding universal
vertex of $G$ extends the preceding small-order immersion.
For an edge $xy$ of $F$, triangle-freeness gives
\begin{equation}\label{eq:small-common}
  \abs{N_G(x)\cap N_G(y)}=n-d_F(x)-d_F(y).
\end{equation}

\emph{Seven vertices.}
The other five vertices supply a totally odd strong $K_3$-immersion.
Lemma~\ref{lem:odd-pair} and~\eqref{eq:small-common} force
$d_F(x)+d_F(y)\ge6$ on every edge. Since $F$ has no isolated
vertex and $\Delta(F)\le3$, it is cubic on seven vertices,
which is impossible.

\emph{Nine vertices.}
The other seven vertices supply a totally odd strong $K_4$-immersion.
Thus $d_F(x)+d_F(y)\ge7$ on every edge.
All degrees are three or four, and the degree-three vertices
are independent. Since $\alpha(F)\le4$, the graph $F$ is not
bipartite. Equation~\eqref{eq:cycle-degree} and $\delta(F)\ge3$
show that a shortest odd cycle is a five-cycle
$v_0v_1v_2v_3v_4v_0$, with indices modulo five.
At most two cycle vertices have degree three, so
\[
  18\le\sum_{i=0}^4d_F(v_i)\le18.
\]
Equality means that exactly two nonadjacent cycle vertices have
degree three, and all four outside vertices have two cycle neighbours.
Let $A_i$ consist of $v_i$ and the vertices whose cycle
neighbourhood is $\{v_{i-1},v_{i+1}\}$.
Triangle-freeness permits edges only between consecutive parts;
maximality makes consecutive parts complete to one another,
since a missing such edge would have no possible common neighbour.
Writing $a_i=\abs{A_i}$, the equations
$a_{i-1}+a_{i+1}=d_F(v_i)$ give, after rotation and reflection,
\[
  (a_0,a_1,a_2,a_3,a_4)=(2,1,2,2,2).
\]
Write $A_0=\{x,x'\}$, $A_1=\{y\}$ and $A_4=\{z,z'\}$,
and choose $w\in A_2$. In $G$, use branch set
$A_0\cup A_1\cup A_3$.
Its only missing pairs are $xy,x'y$; route them by
$x,w,z,y$ and $x',w,z',y$, and use direct edges otherwise.
This is a totally odd strong $K_5$-immersion, a contradiction.

\emph{Twelve vertices.}
If a nonadjacent pair $x,y$ in $G$ has at least three common
neighbours, omit $x,y$ and one common neighbour $z$.
The remaining nine vertices have a totally odd strong $K_5$-immersion.
Since $z$ is outside its branch set, the second sufficient
condition of Lemma~\ref{lem:odd-pair} extends it to $K_6$.
Consequently~\eqref{eq:small-common} forces
$d_F(x)+d_F(y)\ge10$ on every edge of $F$.
Since $\Delta(F)\le5$ and there are no isolated vertices, $F$
is $5$-regular. It is nonbipartite because $\alpha(F)\le5$.
But a shortest odd cycle has length at least five, contradicting
\eqref{eq:cycle-degree}, which would give $25\le24$.
\end{proof}

\subsection{Proof of the theorem}
\label{sec:appendix-induction}

Every vertex outside an induced five-cycle is adjacent to three
consecutive cycle vertices. Indeed, if it misses $v_j$, it must
see $v_{j+2},v_{j+3}$ and at least one of $v_{j+1},v_{j+4}$;
otherwise an independent triple results.

\begin{lemma}
\label{lem:odd-cycle}
Let $C=v_0v_1v_2v_3v_4v_0$ be an induced five-cycle.
Suppose $G-V(C)$ has a totally odd strong $K_m$-immersion with
branch paths of length at most three and branch set $M$, and put
$W=V(G)\setminus(V(C)\cup M)$.
If
\[
  \abs W\ge m,\qquad d_G(v_i)\ge m+2\quad(0\le i<5),
\]
then the immersion extends to a totally odd strong $K_{m+2}$-immersion
with branch paths of length at most three, retaining all old paths.
\end{lemma}

\begin{proof}
All indices are modulo five. Put $w=\abs W$, and define
\[
\begin{aligned}
  P_i&=\{u\in M:N(u)\cap V(C)=\{v_{i-1},v_i,v_{i+1}\}\},\\
  Z_i&=\{y\in W:N(y)\cap V(C)=\{v_{i-1},v_i,v_{i+1}\}\}.
\end{aligned}
\]
These classes contain only vertices with exactly three cycle
neighbours; vertices with four or five neighbours are left
unclassified. Write $p_i=\abs{P_i}$, $z_i=\abs{Z_i}$, and
\[
  A_j=M\setminus N(v_j),\qquad Q_j=W\cap N(v_j),\qquad q_j=\abs{Q_j}.
\]
Each cycle vertex has exactly two neighbours on $C$, so the
degree hypothesis gives
\begin{equation}\label{eq:cycle-counts}
  q_j\ge\abs{A_j},
  \qquad
  \sum_i p_i\le m\le w,\quad
  \sum_i z_i\le w,\quad
  \sum_i q_i\ge3w.
\end{equation}
The last inequality holds because every outside vertex sees
at least three cycle vertices.

Try $a=v_{i-1},b=v_{i+1}$ as the new roots.
A missing pair $a,u$ can use any connector $y\in Q_{i-1}$,
except when $u\in P_{i+2}$ and $y\in Z_i$.
Indeed, the required path is $a,y,v_j,u$ with
$v_j\in(N(y)\cap N(u)\cap V(C))\setminus\{a,b\}$.
If $u$ has four cycle neighbours, it sees all three nonroot
cycle vertices, so such a $v_j$ exists.
Otherwise $u$ has type $P_{i+1}$ or $P_{i+2}$.
In the first case every connector seeing $a$ works; in the
second the unique incompatible cycle neighbourhood of $y$ is
$\{v_{i-1},v_i,v_{i+1}\}$.
The statement at $b$ is the reflected one.
For each $u\in A_{i-1}$, choose a vertex of $Q_{i-1}$, excluding
$Z_i$ when $u\in P_{i+2}$.
Similarly, for each $u\in A_{i+1}$, choose a vertex of $Q_{i+1}$,
excluding $Z_i$ when $u\in P_{i-2}$.
All chosen vertices must be distinct.

We apply Hall's theorem.
At either root, the sets from which we choose are nested, and
$q_j\ge\abs{A_j}$ by \eqref{eq:cycle-counts}.
For choices involving both roots, observe that
\[
  Q_{i-1}\cup Q_{i+1}=W,\qquad
  Z_i\subseteq Q_{i-1}\cap Q_{i+1}.
\]
Their union is therefore $W$, unless both exclude $Z_i$, when
it is $W\setminus Z_i$.
Since $A_{i-1}$ and $A_{i+1}$ are disjoint and $m\le w$,
Hall's condition can fail only if at least one of
\begin{align}
  p_{i+2}&>q_{i-1}-z_i, &&\text{denoted }L_i,\label{eq:failure-L}\\
  p_{i-2}&>q_{i+1}-z_i, &&\text{denoted }R_i,\label{eq:failure-R}\\
  p_{i+2}+p_{i-2}&>w-z_i, &&\text{denoted }T_i\label{eq:failure-T}
\end{align}
holds.

The sets $Z_i,Z_{i-2}$ are disjoint subsets of $Q_{i-1}$,
and $P_{i+1},P_{i+2}$ are disjoint subsets of $A_{i-1}$.
Hence $L_i$ and $R_{i-2}$ cannot both hold: adding them gives
\[
  p_{i+1}+p_{i+2}
    >2q_{i-1}-z_i-z_{i-2}
    \ge q_{i-1}\ge\abs{A_{i-1}}.
\]
Nor can $L_i$ and $T_{i-2}$ both hold, since their sum gives
\[
  p_i+p_{i+1}+p_{i+2}
    >q_{i-1}+w-z_i-z_{i-2}\ge w\ge m.
\]
By reflection, $R_i$ excludes $L_{i+2}$ and $T_{i+2}$.

Suppose every root pair fails, and select one of its
obstructions. If any selected obstruction is $L_i$, the one
at $i-2$ must also be an $L$; iterating forces all five to
be $L$'s. Similarly, an $R$ forces all five to be $R$'s.
The only other possibility is five $T$'s.
Summing the five $L$ inequalities (or the five $R$ inequalities)
gives $\sum_i p_i>2w$, while summing the five $T$ inequalities
gives $2\sum_i p_i>4w$.
Both contradict \eqref{eq:cycle-counts}.

Some root pair therefore has a matching of distinct connectors
to its missing branch pairs. Use the corresponding length-three
paths, all available direct root--branch edges, and the
length-three arc of $C$ between the roots.
Each old branch misses at most one root, and the connectors
are distinct, so the new paths are edge-disjoint. Their
interiors avoid the new branch set, and every new edge meets
$C$, preserving the old immersion.
\end{proof}

\begin{proof}[Proof of Theorem~\ref{thm:two-fifths}]
Induct on $n$, putting $k=\ceil{(2n+2)/5}$.
An edge suffices for $n\le4$.
For $5\le n\le12$, apply Lemma~\ref{lem:odd-small} to an induced
subgraph of order five, seven, nine or twelve as appropriate,
discarding surplus branch vertices if necessary.
Suppose $n\ge13$, and pass to an edge-minimal spanning subgraph
with independence number at most two.
If it has a $K_k$, we are done. Otherwise its nonneighbourhoods
are cliques of size at most $k-1$, so $\delta(G)\ge n-k$.
The graph is connected: a disconnected graph with $\alpha\le2$
consists of at most two cliques, one of order
at least $\ceil{n/2}\ge k$.

There is an induced five-cycle. To see this, take an induced
path $u,v,w$, which exists since $G$ is connected and noncomplete.
Edge-minimality supplies a vertex $s$ nonadjacent to both $u,v$
and a vertex $t$ nonadjacent to both $v,w$.
They are distinct, since otherwise $u,w,s$ would be independent.
The condition $\alpha(G)\le2$ forces $ws,ut,st$ to be edges,
giving the induced cycle $C=u,v,w,s,t,u$.

The graph $G-V(C)$ has at least eight vertices and hence an edge.
Induction gives a totally odd strong $K_m$-immersion there,
where $m=\ceil{(2(n-5)+2)/5}=k-2$.
There are $w=n-k-3$ other vertices outside $C$.
For $n\ge13$ we have $k\le\floor{(n-1)/2}$, so $w\ge m$.
Moreover,
\[
  d_G(v_i)\ge n-k=w+3\ge m+2\qquad(v_i\in V(C)).
\]
Lemma~\ref{lem:odd-cycle} extends the immersion to $K_{m+2}=K_k$.
\end{proof}

\end{document}